\documentclass[12pt, titlepage]{article}
\usepackage{amssymb,amsmath,amsfonts,amssymb}
\usepackage[top=16mm, bottom=10mm, left=35mm, right=28mm]{geometry}
\newcommand{\NI}{\noindent}
\newtheorem{theorem}{\NI{\bf Theorem}}[section]
\newtheorem{lemma}{\NI\bf Lemma}[section]

\newtheorem{cor}{\NI\bf Corollary}[section]
\newtheorem{Rem}{\NI\bf Remark}[section]

\newtheorem{proof}{\NI \bf Proof}
\newcommand{\bp}{\begin{proof}}
\newcommand{\ep}{\end{proof}}
\newcommand{\bt}{\begin{theorem}}
\newcommand{\et}{\end{theorem}}
\newcommand{\bc}{\begin{cor}}
\newcommand{\ec}{\end{cor}}
\newcommand{\bea}{\begin{eqnarray}}
\newcommand{\eea}{\end{eqnarray}}
\newcommand{\ben}{\begin{eqnarray*}}
\newcommand{\een}{\end{eqnarray*}}

\newcommand{\vsp}{\vskip 1em}

\def \qed {\hfill \vrule height 6pt width 6pt depth 0pt}
\newcommand{\be}{\begin{equation}}
\newcommand{\ee}{\end{equation}}

\title {\bf On testing  More IFRA Ordering-II }%Individual Risk Models }
\author{
Muhyiddin Izadi\thanks {Corresponding author}
 \ \ \ \ \  Baha-Eldin Khaledi 
\\Department of Statistics, Razi
University, Kermanshah, Iran
 \and Chin-Diew Lai \\ Institute of Fundamental Sciences - Statistics and Bioinformatics,\\ Massey University, Palmerston North, New Zealand }
    \date{}
\begin {document}
\maketitle

\begin{abstract}
Suppose $F$ and $G$ are two life distribution functions. It is
said that $F$ is more IFRA than $G$ (written by $F\le_*G$) if
$G^{-1}F(x)$ is starshaped on $(0, \infty)$. In this paper, the problem
of testing $H_0: F=_*G$ against $H_1: F\le_*G$ and $F\neq_*G$ is
considered in both cases when $G$ is known and when $G$ is
unknown. We propose a new test based on U-statistics and obtain the asymptotic distribution of the test
statistics. The new test is compared with some well known tests in the literature. In addition, we apply our test to a real data set in the context of reliability.

 \vsp {{\bf Keywords} : Asymptotic normality,  star order, increasing failure rate average, Pitman's asymptotic efficiency,  U-statistic.}

\end{abstract}

%{\bf Keywords : }.
 \maketitle
% -----------------------------------------------------------
\section{Introduction}
Let $X$ be a lifetime of an appliance with density function $f$,
distribution function $F$ and survival function $\bar F$. Let also denote $F^{-1}$ as the right continuous inverse function of $F$.  $X$ is said to be IFRA ({\it increasing failure rate average}) if $\tilde{r}_F(x)=\frac{\int_{0}^{x}r_F(t)dt}{x}$
is nondecreasing in $x\ge0$ which is equivalent to that $-\frac{\log\bar{F}(x)}{x}$ is nondecreasing in $x\ge0$ where
 $r_F(x)=\frac{f(x)}{\bar{F}(x)}$. It is of considerable interest to producers and
users of the appliances to evaluate the severity of average failure risk 
at a particular point of time and to see if $\tilde{r}_F(x)$ is either
increasing or decreasing in time. That is, it is of practical
importance to characterize the aging class of underlying random
lifetimes. In particular, since the IFRA class of aging is one of
the most important aging classes, testing that the
distribution $F$ has a constant hazard rate against the hypothesis
that $F$ is IFRA has been studied extensively in the literature;
see for example,  Deshpande (1983), Kochar (1985), Link (1989),
Ahmad (2000) and El-Bassiouny (2003) among others.
% Now we introduce some of the
%know tests in the literature that we use to compare with the new
%proposed test in this paper.
In fact, $F$ is IFRA if and only if $\frac{E_\lambda^{-1}F(x)}{x}$ is nondecreasing in $x\ge0$ or equivalently $\frac{\tilde{r}_F(F^{-1}(u))}{\tilde{r}_E(E_\lambda^{-1}(u))}$ is nondecreasing in $u\in (0,1)$ where $E_\lambda$ is an exponential distribution with mean $\lambda$. This implies that $F$ ages faster than $E$, i.e., $F$ is more IFRA than $E_\lambda$.

In order to evaluate the performance of an appliance, we need to
compare its aging behavior with some distributions other
than exponential distribution such as the Weibull,
gamma,  linear failure rate  or even an unknown distribution $G$. The notion of the {\it star order} that establishes an equivalent
class of distributions is one of the useful tools for this
comparison. Let $Y$ be another non-negative random variable
with distribution function $G$. We
say that $X$ is less than $Y$ with respect to the {\it star order} (written by $X\le_*Y$ or $F\le_*G$ ) if
$G^{-1}F(x)$ is starshaped on $[0,\infty)$; that is, $\frac{G^{-1}F(x)}{x}$ is nondecreasing in $x\ge0$. It is known that
\be
\label{*-order} F\le_*G \Leftrightarrow
\frac{\tilde{r}_{F}(F^{-1}(u))}{\tilde{r}_G(G^{-1}(u))} \ \ \mbox{is nondecreasing
in}\ u\in(0,1),
 \ee
where $\tilde{r}_F$ and $\tilde{r}_G$ are failure rate average functions of
$F$ and $G$, respectively. Using $(\ref{*-order})$, the relation
$X\le_*Y$ is interpreted as $X$ ages faster than $Y$ and it is
said that $X$ is more IFRA than $Y$  (cf. Kochar and Xu, 2011 ). It is obvious that if $F\le_*G$ and $G\le_*F$ then $F(x)=G(ax)$ for all $x\ge0$ and some $a>0$. In this case, we say $F=_*G$.

Izadi and Khaledi (2012) have considered the problem of testing the null hypothesis $H_0: F=_*G $ against $H_1: F\le_*G $ and $ F\neq_*G$. They proposed a test based on kernel density estimation.
In this paper, we further study this problem of testing in the one-sample as well as the two-sample problem and propose a new simple test based on a U-statistic. In  both cases, we compare the new proposed test with some well known tests in the literature. It is found that  our test is comparable to  the others.

To establish our new test we need the following lemma.

\begin{lemma}\label{StarLemma}
Let $X_1, X_2$ $(Y_1, Y_2)$ be two independent copies of the random variable $X$ ($Y$) with
distribution function $F$ ($G$) and let $\mu^{(2)}_F=E[
\max\{X_1,X_2\}]$ ($\mu^{(2)}_G=E[\max\{Y_1,Y_2\}]$) where $E[.]$ is the expectation operator. If $F$ is more IFRA than $G$, then
\begin{eqnarray*}\label{inequality1}
   \frac{\mu^{(2)}_F}{\mu_F}\le \frac{\mu^{(2)}_G}{\mu_G}
 \end{eqnarray*}
where $\mu_F$ ($\mu_G$) is the expectation of $F$ ($G$).
\end{lemma}
{\bf Proof:} We know that more IFRA order is scale invariant. Thus, $X\le_*Y$ implies $\frac{\mu_Y}{\mu_X}X\le_*Y$. Now, the required result follows from Theorem 7.6 of Barlow and Proschan (1981, page 122).

\begin{Rem}
The above lemma has been proved by  Xie and Lai (1996) under the condition that $F$ is more IFR than $G$ (for definition, see Shaked and Shantikumar, 2007, p. 214) which is stronger than more IFRA order.
\end{Rem}

 Now, let $\delta_F=\frac{\mu^{(2)}_F}{\mu_F}$, $\delta_G= \frac{\mu^{(2)}_G}{\mu_G}$ and
\begin{eqnarray}\label{measure}
%   \delta(F, G)
    \delta(F, G)=\delta_F-\delta_G.
 \end{eqnarray}
 It is obvious that if $F=_*G$, then $\delta(F,G)=0$ and  if $F\le_*G$ and
$F\neq_*G$, then it follows  from Lemma \ref{StarLemma} that $\delta(F,G)<0$. That is, $\delta(F,G)$   can be
considered as a measure of departure from $H_0:
 F=_*G$ in favor of $H_1: F\le_*G$ and $F\neq_*G$. So, our test statistic is based on the estimation of $\delta(F,G)$.

The organization of this paper is as follows.  In Section 2, we
propose the new test for the case when $G$ is known.  The case
when $G$ is unknown is studied in Section 3. In Section 4, the
performance of our test is evaluated and compared.

\section{The One-Sample Problem}
%====================================================================== New Section Test statistic

Let  $G_0$ be a known  distribution function and $X_1, \ldots,
X_n$ be a random sample from an unknown distribution $F$. Now by
using the measure (\ref{measure}), the test statistic
\begin{eqnarray*}
  \hat\delta(F, G_0)&=&\hat\delta_F- \delta_{G_0}\nonumber
\end{eqnarray*}
is used for testing
\begin{eqnarray*}
H_0: F=_*G_0
\end{eqnarray*}
against
\begin{eqnarray*}
H_1: F\le_*G_0 \ \mbox{and} \ F\neq_*G_0
\end{eqnarray*}
 where
\begin{eqnarray}\label{delta_F}
\hat\delta_F&=&\frac{\displaystyle{\mathop{\sum\sum}_{i\neq j}}\max\{X_i, X_j\}}{n(n-1)\bar X} 
\end{eqnarray}
and $\bar{X}$ is the mean of the random sample. %Note that in (\ref{delta_F}), $\displaystyle\sum$ denotes a double summation.
In the next theorem, we obtain the asymptotic distribution of
$\hat\delta(F, G_0)$ by using the standard theory of
U-statistics.
\begin{theorem}\label{theorem}
Suppose  $E[\max\{X_1, X_2\}-\frac{\delta_F}{2}(X_1+X_2)]^2<\infty$. As $n\rightarrow\infty$, $n^{1/2}[\hat\delta(F, G_0)-\delta(F,
G_0)]$ is asymptotically normal with mean $0$ and variance
\begin{eqnarray}\label{Sigma}
    \sigma_F^2=\frac4{\mu^2_F}\times
    Var\left(XF(X)+\int_{X}^{\infty}tdF(t)-\frac{\delta_F}{2}X\right).
\end{eqnarray}
Under $H_0$, $n^{1/2}\hat\delta(F, G_0)$ is asymptotically
normal with mean $0$ and variance $\sigma^2_0=\sigma^2_{G_0}.$
\end{theorem}
%========================================================= Proof
{\bf Proof:} First note that
\begin{eqnarray*}
  \hat\delta_F-\delta_F &=& \frac{\displaystyle{\mathop{\sum\sum}_{i\neq j}}\left[\max\{X_i,
  X_j\}-\frac{\delta_F}{2}(X_i+X_j)\right]}{n(n-1)\bar X}\\
  &=&\frac{\displaystyle{\mathop{\sum\sum}_{i\neq j}}\phi(X_i, X_j)}{n(n-1)\bar X},
\end{eqnarray*}
where
\begin{eqnarray*}
\phi(X_i, X_j)=\max\{X_i,
  X_j\}-\frac{\delta_F}{2}(X_i+X_j).
\end{eqnarray*}
Let define
\begin{eqnarray*}
T^*=\frac{\displaystyle{\mathop{\sum\sum}_{i\neq j}}\phi(X_i, X_j)}{n(n-1)}.
\end{eqnarray*}
 By the standard theory
of U-statistics, if $E[\phi^2(X_1,X_2)]<\infty$, as $n\longrightarrow\infty$
\begin{eqnarray*}
  \frac{\sqrt nT^*_n}{\sigma_*}\stackrel d
  \rightarrow N(0,1)
\end{eqnarray*}
where
\begin{eqnarray*}\label{variance}
  \sigma_*^2=4\times Var(\phi_1(X))
\end{eqnarray*}
and
\begin{eqnarray*}
\phi_1(x)&=&E[\phi(x, X)].
\end{eqnarray*}
Now by the strong law of large numbers we have $\bar X\stackrel
{a.s.}\rightarrow\mu_F$ and hence, by Slutsky theorem $\sqrt
n[\hat\delta(F, G_0)-\delta(F, G_0)]$ is asymptotically normal
with mean 0 and variance $\sigma^2_F=\frac{\sigma^2_*}{\mu^2_F}$.
Under $H_0$, $\delta(F, G_0)=0$ and $\sigma^2_0=\sigma^2_{G_0}$. \
\ \ \ $\square$ \vspace{0.6cm}
%===================================================================== End of Proof
%=============================================================================

A small value of $\hat\delta(F, G_0)$ indicates that testing
$H_0$ against $H_1$ is significant. Thus, we reject $H_0$ at level $\alpha$ if
$n^{1/2}\hat\delta(F, G_0)/\sigma_{G_0}<z_{\alpha}$, where
$z_{\alpha}$ is $\alpha^{th}$ quantile of the standard normal distribution.

In the case  $G_0(x)=E_\lambda(x)=1-\exp(-\lambda x)$, $x\ge0$ and $\lambda>0$, the problem  is testing the null hypothesis  $H_0: F$ is an exponential distribution
against the alternative hypothesis $H_1: F$ is IFRA and not exponential.
It can be shown that $\delta_{E_\lambda}=\frac32$ and $\sigma^2_{E_\lambda}=\frac1{12}$.
By the above theorem, under $H_0$, $\sqrt n(\hat\delta_{F}-3/2)$ is asymptotically normal with mean 0 and variance $\frac1{12}.$
Thus we reject $H_0$ in favor of $H_1$ if $\sqrt {12n}(\hat\delta_{F}-3/2)<
 z_{\alpha}$. 

 In the following we find the exact distribution of $\hat{\delta}_F$ under the hypothesis $F$ is an exponential distribution. First, note that we can rewrite $\hat{\delta}_F$ as
 \begin{eqnarray}
 \hat{\delta}_F=\frac{2\displaystyle\sum_{i=1}^{n}(i-1)X_{(i)}}{n(n-1)\bar{X}}= \frac{\displaystyle\sum_{i=1}^{n}c_{i:n}D_i}{\displaystyle\sum_{i=1}^{n}D_i}
 \end{eqnarray}
 where $X_{(i)}$ is the $i^{th}$ order statistic of $X_i$'s, $$D_i=(n-i+1)(X_{(i)}-X_{(i-1)}), \ \  \ c_{i:n}=\frac{2\displaystyle\sum_{j=i}^{n}(j-1)}{(n-1)(n-i+1)}$$
 and assuming $X_{(0)}=0$. Now, by the same arguments as in  Langenberg and Srinivasan (1979), we will get the following result.
 \begin{theorem}\label{Exactdist}
 Let $F$ be an exponential distribution, then 
 \begin{eqnarray}
 P\{\hat{\delta}_F\le x\}=1-\sum_{i=1}^{n}\mathop{\prod_{j=1}^{n}}_{j\neq i}\frac{c_{i:n}-x}{c_{i:n}-c_{j:n}}I{(x<c_{i:n})}
 \end{eqnarray}
 where $I(.)$ is the usual indicator function.
 \end{theorem}

\begin{table}[!ht] 
\centering 
\small 
\caption{ {Critical values of $\sqrt{12n}(\hat{\delta}_F-3/2)$ for small sample sizes}} \label{TabCritic} 
\vspace{0.66cm}
\small\addtolength{\tabcolsep}{-3pt}

\begin{tabular}{cccc|ccc}
%\hrule
%\hline 
n & & $\alpha$: Lower Tail & & &$\alpha$: Upper Tail  &\\ 
  & 0.01 &0.05&0.1 &0.1&0.05&0.01\\
  \hline
2 &-2.400500 & -2.204541 &-1.959592 & 1.959592 & 2.204541 & 2.400500\\
3&-2.575752 &-2.051328 &-1.658360 & 1.658360 & 2.051328 & 2.575752\\
4&-2.560006 &-1.918143 &-1.516280  &1.516280  &1.918143  &2.560006\\
5&-2.517587 &-1.846175 &-1.458997  &1.458997  &1.846175  &2.517587\\
6&-2.482569 &-1.807959 &-1.424446  &1.424446  &1.807959  &2.482569\\
7& -2.458901 &-1.781575 &-1.400767  &1.400767  &1.781575  &2.458901\\
8&-2.441786 &-1.762473 &-1.383974  &1.383974  &1.762473  &2.441786\\
9&-2.428500 &-1.748106 &-1.371312  &1.371312  &1.748106  &2.428500\\
10&-2.417939 &-1.736865 &-1.361442 & 1.361442 & 1.736865 & 2.417939\\
11&-2.409356 &-1.727862 &-1.353531  &1.353531  &1.727862 & 2.409356\\
12&-2.402239 &-1.720450 &-1.347047  &1.347047  &1.720450  &2.402239\\
13&-2.396243 &-1.714193 &-1.341635  &1.341635  &1.714193  &2.396243\\
14&-2.391124 &-1.708937 &-1.337050  &1.337050  &1.708937  &2.391124\\
15&-2.386703 &-1.704422 &-1.333116  &1.333116  &1.704422  &2.386703\\
16&-2.382846 &-1.700502 &-1.329703  &1.329703  &1.700502  &2.382846\\
17&-2.379451 &-1.697066 &-1.326714  &1.326714  &1.697066  &2.379451\\
18&-2.376441 &-1.694029 &-1.324074  &1.324074  &1.694029  &2.376441\\
19&-2.373754 &-1.691327 &-1.321727  &1.321727  &1.691327  &2.373754\\
20&-2.371340 &-1.688906 &-1.319625  &1.319625  &1.688906  &2.371340\\
21&-2.369160 &-1.686725 &-1.317732  &1.317732  &1.686725  &2.369160\\
22&-2.367182 &-1.684749 &-1.316018  &1.316018  &1.684749  &2.367182\\
23&-2.365378 &-1.682952 &-1.314460  &1.314460  &1.682952  &2.365378\\
24&-2.363726 &-1.681309 &-1.313036  &1.313036  &1.681309  &2.363726\\
25&-2.362209 &-1.679803 &-1.311730  &1.311730  &1.679803  &2.362209\\
26&-2.360810 &-1.678415 &-1.310529  &1.310529  &1.678415  &2.360810\\
27&-2.359516 &-1.677134 &-1.309419  &1.309419  &1.677134  &2.359516\\
28&-2.358316 &-1.675947 &-1.308392  &1.308392  &1.675947  &2.358316\\
29&-2.357199 &-1.674844 &-1.307437  &1.307437  &1.674844  &2.357199\\
30&-2.356158 &-1.673817 &-1.306548  &1.306548  &1.673817  &2.356158\\
31&-2.355185 &-1.672857 &-1.305718  &1.305718  &1.672857  &2.355185\\
32&-2.354273 &-1.671960 &-1.304942  &1.304942  &1.671960  &2.354273\\
33&-2.353418 &-1.671118 &-1.304214  &1.304214  &1.671118  &2.353417\\
34&-2.352613 &-1.670326 &-1.303529  &1.303529  &1.670326  &2.352612\\
35&-2.351855 &-1.669581 &-1.302885  &1.302885  &1.669581  &2.351854\\
36&-2.351140 &-1.668878 &-1.302278  &1.302278  &1.668878  &2.351138\\
37&-2.350454 &-1.668214 &-1.301704  &1.301704  &1.668214  &2.350461\\
38&-2.349821 &-1.667586 &-1.301161  &1.301161  &1.667586  &2.349821\\
39&-2.349184 &-1.666990 &-1.300647  &1.300647  &1.666991  &2.349213\\
40&-2.348688 &-1.666426 &-1.300159  &1.300159  &1.666426  &2.348636\\
\end{tabular}
%\end{eqnarray}
\end{table}

 By using Theorem \ref{Exactdist}, we tabulate the critical point of $\sqrt{12n}(\hat{\delta}_F-3/2)$ under  exponentiality  for small sample sizes ($\le 40$) in Table \ref{TabCritic}. So, for small sample sizes, we reject exponentiality in favor of IFRA-ness if  $\sqrt{12n}(\hat{\delta}_F-3/2)$ is smaller than the critical point in Table \ref{TabCritic} corresponding with the level of significance chosen.

 El-Bassiouny (2003) has considered the problem of testing exponentiality against IFRA-ness in the alternative and  proposed a class of test. His test is based on the test statistics
 \begin{eqnarray*}
 \hat{\Delta}_{r+1}&=&\frac{2\displaystyle{\mathop{\sum\sum}_{i<j}} \left(\min\{X^{r+1}_i, X^{r+1}_j\}-\frac{X^{r+1}_i}{2}\right)}{n(n-1){\bar X}^{r+1}}
 \end{eqnarray*} 
 and  large values of $\hat{\Delta}_{r+1}$ are significant for the considered problem of testing. 
 If $r=0$,
 \begin{eqnarray}
 \hat{\Delta}_1&=&\frac{2\displaystyle{\mathop{\sum\sum}_{i<j}} \left(\min\{X_i, X_j\}-\frac{X_i}{2}\right)}{n(n-1)\bar X}\nonumber\\
 &=& \frac{2\displaystyle{\mathop{\sum\sum}_{i<j}}\min\{X_i, X_j\}}{n(n-1)\bar X}-\frac12
 \end{eqnarray} 
 On the other hand, using the fact that $\delta_{E_\lambda}=\frac{3}{2}$,
 %In this case, our  test is the same as the test of El-Bassiouny (2003) based on $\hat{\Delta}_1$, because
\begin{eqnarray*}
 \hat\delta(F,E_\lambda)&=&\frac{\displaystyle{\mathop{\sum\sum}_{i\neq j}}\max\{X_i,
  X_j\}}{n(n-1)\bar
  X}-\delta_{E_\lambda}\\
  &=& \frac{2\displaystyle{\mathop{\sum\sum}_{i< j}}\left[(X_i+X_j)-\min\{X_i,
  X_j\}\right]}{n(n-1)\bar
  X}-\frac{3}{2}\\
  &=&\frac{2\left[\displaystyle{\mathop{\sum\sum}_{i<j}}(X_i+X_j)-\mathop{\sum\sum}_{i<j}\min\{X_i,
  X_j\}\right]}{n(n-1)\bar
  X}-\frac{3}{2}\\
  &=&\frac{2\left[(n-1)\displaystyle\sum_{i=1}^{n}X_i-\mathop{\sum\sum}_{i<j}\min\{X_i,
  X_j\}\right]}{n(n-1)\bar
  X}-\frac{3}{2}\\
  &=&\frac{1}{2}-\frac{2\displaystyle{\mathop{\sum\sum}_{i<j}}\min\{X_i,
  X_j\}}{n(n-1)\bar
  X}\\
  &=&-\hat{\Delta}_1.
 \end{eqnarray*}
That is, for the case when $r=0$ and $G_0$ is an exponential distribution, the proposed test is equivalent to that of El-bassiouny (2003).

It is worth to mention that our test is consistent; that is, if $\beta_n(F)$ is the power of our test, then under the alternative hypothesis, $\lim_{n\rightarrow\infty}\beta_n(F)=1$
which follows from Theorem \ref{theorem} and Problem 2.3.16 in Lehmann (1999),
%=================================================== Section Two sample
\section{ The Two-Sample Problem}
In this section, we consider the two-sample problem when $G$ is unknown. Let $X_1, X_2, \ldots, X_n$ and $Y_1, Y_2, \ldots, Y_m$ be two
independent random samples from unknown distribution functions $F$
and $G$, respectively, and $N=n+m$. Assume that $\hat\delta_F$ is as in (\ref{delta_F}) and $\hat\delta_G$
is defined similarly in terms of $Y_1, \ldots, Y_m$ and $\bar
Y$.
%We assume that for all $n$ and
%$m$ the inequalities
%\begin{eqnarray}\label{inequality}
  %0<c\leq\frac{n}{N}\leq1-c<1, \ \ 0<c\le1/2
%\end{eqnarray} hold true.  
The test statistic
\begin{eqnarray*}
\hat\delta(F, G)= \hat\delta_F-\hat\delta_G
\end{eqnarray*}
which is the estimate of the measure in (\ref{measure}) is used for testing the null hypothesis
\begin{eqnarray}\label{H_0G}
  H_0: F=_*G
\end{eqnarray}
against the alternative hypothesis
\begin{eqnarray}\label{H_1G}
  H_1: F\le_*G \ \  \mbox{and} \ \ F\neq_*G.
\end{eqnarray}
Small values of  $\hat\delta(F, G)$ are significant for testing $H_0$ against $H_1$.
In the following theorem we obtain the asymptotic distribution of
$\hat\delta(F, G)$. 
%====================================================================== Theorem
\begin{theorem} If $E[\max\{X_1, X_2\}-\frac{\delta_F}{2}(X_1+X_2)]^2$ and 
$E[\max\{Y_1, Y_2\}-\frac{\delta_G}{2}(Y_1+Y_2)]^2$ are finite and $n$ and $m$ $\rightarrow \infty$  such  that 
$\frac{n}{N}$ $\rightarrow c$, $c\in(0,\frac{1}{2}]$, then $\sqrt N (\hat\delta(F, G)-\delta(F,G))$ 
is asymptotically normal with mean 0 and variance
  \begin{eqnarray*}
    \sigma^2_{F,G}=\frac Nn{\sigma^2_F} + \frac
    Nm{\sigma^2_G}
  \end{eqnarray*}
  where $\sigma^2_F$ is given in (\ref{Sigma}) and $\sigma^2_G$ is
  defined similarly in terms of $Y$.
\end{theorem}
%====================================================================== Proof
%=========================================================================
{\it\bf Proof:} It is easy to see that
\begin{eqnarray}
  \frac{\sqrt N (\hat\delta(F, G)-\delta(F,
  G))}{\sigma_{F,G}}&=&\frac{\sqrt
  m\sigma_F}{\sqrt{m\sigma^2_F+n\sigma^2_G}}[\sqrt
  n(\hat\delta_F-\delta_F)/\sigma_F]\nonumber\\
  &&-\frac{\sqrt
  n\sigma_G}{\sqrt{m\sigma^2_F+n\sigma^2_G}}[\sqrt
  m(\hat\delta_G-\delta_G)/\sigma_G]\nonumber.
\end{eqnarray}

 From  the result of Theorem 2.1, as both $n$ and $m$ $\rightarrow\infty$, we have that
\begin{eqnarray}
\sqrt n(\hat\delta_F-\delta_F)/\sigma_F\stackrel d\rightarrow
N(0,1) \mbox{ and }
\sqrt m(\hat\delta_G-\delta_G)/\sigma_G\stackrel d\rightarrow
N(0,1).\nonumber
\end{eqnarray}
Since $\hat\delta_F$ and $\hat\delta_G$ are independent, the required result follows from the fact that convergence in distribution is closed under the convolution of independent sequences of random variables (cf. Theorem 6.6 of Gut (2009), page 169).   \qed
 %================================================================ End of Proof
 %========================================================================

In practice $\sigma^2_{F,G}$ is unknown, but it can be estimated
by the consistent estimator
\begin{eqnarray}\label{SigmaFG}
  \hat\sigma^2_{F,G}=\frac Nn{\hat\sigma^2_F} + \frac
  Nm{\hat\sigma^2_G}
\end{eqnarray}
where
\begin{eqnarray*}
  \hat\sigma^2_F&=&4\times\frac{\displaystyle\sum_{i=1}^{n}{\hat\phi}^2(X_i)}{n{\bar X}^2},\\
  \hat\phi(X_i)&=&\frac1{n-1}{\small\mathop{\sum_{j=1}^{n}}_{j\neq i}} [\max\{X_i, X_j\}-\frac{\hat\delta_F}{2}(X_i+{X_j})]
\end{eqnarray*}
and
\begin{eqnarray*}
  \hat\sigma^2_G&=&4\times\frac{\displaystyle\sum_{i=1}^{m}{\hat\eta}^2(Y_i)}{m{\bar Y}^2},\\
  \hat\eta(Y_i)&=&\frac1{m-1}\mathop{\sum_{j=1}^{m}}_{j\neq i} [\max\{Y_i, Y_j\}-\frac{\hat\delta_G}{2}(Y_i+{Y_j})].
\end{eqnarray*}
Now by  Slutsky theorem, under $H_0$, $\sqrt N\hat\delta(F,
G)/\hat\sigma_{F,G}$ is asymptotically normal with mean 0 and
variance 1 as both $n$ and $m$ $\rightarrow\infty$. Hence, for large sample sizes, 
%an asymptotically distribution free test for testing $H_0$ of
%(\ref{H_0G}) against $H_1$ of (\ref{H_1G}) is to reject 
$H_0$  is rejected at
level $\alpha$ if $\sqrt N\hat\delta(F, G)/\hat\sigma_{F,G}<
z_{\alpha}$.

\section{Simulation Study}
In this section, we study the performance of  our test and compare it with some well known tests in the literature for the one-sample and the two-sample problems. 
\subsection{The One-Sample}
We recall that in the one-sample problem we consider testing $H_0: F=_*G_0$ against $H_1: F\le_*G_0$ and $F\neq_*G_0$ when $G_0$ is a known distribution. %As mentioned, if $G_0=1-\exp\{-x\}, \ x>0$, then the problem is testing exponentiality against IFRA-ness. In this case, we compare our new test with some well known tests which are introduced below:
For the case when $G_0(x)=1-\exp\{-\lambda x\}, \ x>0$, we compare our proposed test with the following well known tests which are in the literature. Note that in this case the problem is testing exponentiality against IFRA-ness.    

%In the following, We compare our new test with some well known tests for the case when  $G_0=1-\exp\{-x\}, \ x>0$,
\begin{description}
%==================================================================== Deshpande's test
\item Deshpande (1983): The test statistics is 
\begin{eqnarray*}\label{Deshpande's test}
J_b=\frac{1}{n(n-1)}\mathop{\sum\sum}_{i\neq j}h_b(X_i, X_j), \ \ \ b\in (0,1)
\end{eqnarray*}
where $h_b(x,y)=1,$ if $x>by$; 0, otherwise. Large values of $J_b$ are used to reject  exponentiality in favor of IFRA-ness. It has been shown that under $H_0$, $n^{1/2}(J_b-(b+1)^{-1})$ is asymptotically normal with mean zero and variance $4\xi$ where
\begin{eqnarray*}
\xi=\frac{1}{4}\{1+\frac{b}{b+2}+\frac{1}{2b+1}+\frac{2(1-b)}{b+1}-\frac{2b}{b^2+b+1}-\frac{4}{(b+1)^2}\}.
\end{eqnarray*}
Deshpande (1983) has recommended $b=0.9$. 
%=========================================================                        Kochar test
\item Kochar (1985):  $H_0$ is rejected  for large values of 
\begin{eqnarray}
T_n=\frac{\sum_{i=1}^{n}J(\frac{i}{n+1}) X_{(i)}}{n\bar{X}}, \ \ \ \ J(u)=2(1-u)[1-\log(1-u)]-1.
\end{eqnarray}
The asymptotic distribution of $(108n/17)^{1/2}T_n$ is the standard normal distribution.
%=================================================================== Link's test
\item Link (1989): Large values of the test statistic
\begin{eqnarray}
\Gamma=\frac{2}{n(n-1)}\mathop{\sum\sum}_{i<j}\frac{X_{(i)}}{X_{(j)}}.
\end{eqnarray}
certify that $F$ is IFRA. For large values of $n$, under $H_0$, the distribution of $\frac{\sqrt{n}(\Gamma-(2\log2 -1))}{\sqrt{0.048225}}$ is approximately standard normal. 
%================================================================== Ahmad 2000
%===========================================================================
\item Ahmad (2000): The test statistic is 
\begin{eqnarray}
\hat\Delta_F=[n(n-1)a_n]^{-1}\mathop{\sum\sum}_{i\neq j}X_ik\left(\frac{X_i-X_j}{a_n}\right),
\end{eqnarray}
where $k$ is a known symmetric density function and $a_n$ is a sequence of positive real numbers such that $na_n\longrightarrow\infty$ and $na^4_n\longrightarrow0$. Under some conditions, $\sqrt{\frac{108n}{5}}(\hat{\Delta}_F-\frac{1}{4})$ is asymptotically normal with mean zero and variance 1, when $F$ is an exponential distribution (Ahmad, 2000). $H_0$ is rejected at level $\alpha$ if $\sqrt{\frac{108n}{5}}(\hat{\Delta}_F-\frac{1}{4})>z_{1-\alpha}$. Ahmad has recommended standard normal density as kernel function and $a_n=n^{-\frac{1}{2}}$.
%=====================================================================
% ==================================================================== Ahmad's test
%\item Ahmad (1994): A class of test has been introduced. Let $0<a\le1$ and $k\ge1$. The test statistics is 
%\begin{eqnarray}
%\hat{J}(F_n; a, k)={(\frac{a}{a+1})}^k-\frac{\displaystyle\sum_{c}I(aX_{i_1}>X_{i_2}+ \ldots +X_{i_{k+1}})}{n \binom{n-1}{k}}
%\end{eqnarray}
%where $\displaystyle\sum_{c}$ extends over all indices $i_j$ such that $1\le i_2<\ldots<i_{k+1}\le n$ and $i_1\neq i_j$, $j=2, \ldots, k+1$.
%Under $H_0$, $\sqrt{n}(\hat{J}(F_n; a, k))$ is asymptotically normal with mean zero and variance as given (2.5) in Ahmad (1994).
%If $k=1$, this test is equivalent to the test of Deshpande (1983). Ahmad recommended $\alpha=0.4$ and $k\le4$ for testing exponentiality against IFRA-ness.
\end{description}

First, we  investigate the accuracy of normal distribution as the limit distribution of the test statistics under $H_0$.  
%To justify the accuracy of the standard normal distribution as the limit distribution of $\sqrt{12n}(\hat{\delta}_F-3/2)$ under $H_0$,
In order to do this, we simulate the size of  the tests for  nominal sizes $\alpha$= 0.01, 0.05, 0.1 and  large sample sizes $n=40 (5)60(10)70$.  In the simulation, 10000 samples are generated from exponential distribution with mean 1. The calculated size is the proportion of 10000 generated samples that resulted in rejection of $H_0$ where the rejection regions have been obtained by the asymptotical distribution of test statistics. The simulated values were tabulated in Table \ref{Tablesize}.  All simulations were done by R package.

From Table \ref{Tablesize}, we find that the tests by Deshpande (1983) and Kochar (1985) are over shoot the nominal sizes for all sample sizes. The simulated sizes of the tests due to Link (1989) and Ahmad (2000) are greater than the nominal sizes but Link's test always dominates Ahmad's test. 
It is clear from the contents of  Table \ref{Tablesize} that  the simulated sizes of our new test are much closer  to the nominal sizes for all sample sizes.      
\begin{table}[htb] 
\centering 
\small 
\caption{ Simulated sizes of our test for different nominal sizes and large sample sizes}\label{Tablesize} 
\vspace{0.5cm}
\small\addtolength{\tabcolsep}{-3pt}

\begin{tabular}{clccc|clccc}
%\hrule
%\hline 
n & &&  nominal size ($\alpha$) & &n & &&  nominal size ($\alpha$) & \\ 
  && 0.01 &0.05&0.1 &&& 0.01 &0.05&0.1\\
  \hline
40& $\hat{\delta}(F,E)$ & 0.0104 &0.0518 & 0.1044&   55 &$\hat{\delta}(F,E)$&0.0108  &0.0500   &0.1003\\
   &$J_{0.9}$                     &0.0637  &0.1243 &0.1709&    & $J_{0.9}$                  &0.0550  &0.1101   &0.1674\\
   &$T_n$                    &0.0396   &0.1815 &0.3157&    & $T_n$                  &0.0346  &0.1565  &0.2753\\
   &$\Gamma$              &0.0181  &0.0612  &0.1110&   & $\Gamma$           &0.0151  &0.0569   & 0.1063\\
   &$\hat{\Delta}_F$      &0.0311 &0.0772 &0.1196&   &$\hat{\Delta}_F$    &0.0324    &0.0783   &0.1237\\
 45&$\hat{\delta}(F,E)$&0.0106  &0.0505  &0.1015&    60 &$\hat{\delta}(F,E)$&0.0101  &0.0517  &0.1042\\
    & $J_{0.9}$                  &0.0661  &0.1163  &0.1769&     & $J_{0.9}$                  &0.0494  &0.1128 &0.1680\\
    & $T_n$                  &0.0380  &0.1704  &0.2938&     & $T_n$                  &0.0349 &0.1548 &0.2666\\
    & $\Gamma$           &0.0167  &0.0584   &0.1089&     & $\Gamma$           &0.0157 &0.0581  &0.1095\\
    &$\hat{\Delta}_F$    &0.0338 &0.0796 &0.1232&     &$\hat{\Delta}_F$    &0.0304 &0.0742 &0.1209\\
50 &$\hat{\delta}(F,E)$& 0.0112 &0.0494   &0.1009&    70 &$\hat{\delta}(F,E)$& 0.0090 &0.0489 &0.1024\\
    & $J_{0.9}$                  &0.0601  &0.1220   &0.1736&     & $J_{0.9}$                  & 0.0469 &0.1048 &0.1587\\
    & $T_n$                  &0.0371 &0.1654   &0.2852&     & $T_n$                  & 0.0303 &0.1410 &0.2534\\
    & $\Gamma$           &0.0166  &0.0563   &0.1055&     & $\Gamma$           & 0.0147 &0.0580 &0.1072\\  
    &$\hat{\Delta}_F$    & 0.0312 &0.0810 &0.1232&     &$\hat{\Delta}_F$   & 0.0310 &0.0783 &0.1240\\     
\end{tabular}
%\end{eqnarray}
\end{table}

In the following, to assess how our proposed  test performs  relatively, we first consider the large sample sizes and  use the measure of  Pitman's asymptotic relative efficiency (PARE) (cf. Nikitin, 1995, Section 1.4). Consider testing $H_0$ that $F$ is an exponential distribution against $H_1$ that $F=F_{\theta_n}$ where $\theta_n=\theta_0 + kn^{-\frac12}$, $k$ is an arbitrary positive constant and $F_{\theta_0}$ is exponential. Then, Pitman's
asymptotic efficiency (PAE) of a test based on statistic $T_n$ is
\begin{eqnarray}\label{PAEf}
  PAE(T_n)=\lim_{n\rightarrow\infty}\frac{\left[\frac{\partial
  E_{\theta}(T_n)}{\partial\theta}|_{\theta=\theta_0}\right]^2}{Var_{\theta_0}[\sqrt{n}T_n]}.
\end{eqnarray}
 Using (\ref{PAEf}), the PAE of our test is given by

\begin{eqnarray*}
  PAE(\hat\delta(F, E_{\lambda}))&=&\frac{(\frac{\partial
  \delta_{F_{\theta}}}{\partial\theta}|_{\theta=\theta_0})^2}{\sigma^2_{F_{\theta_0}}}.
%  &=& 12\left[\int_{0}^{\infty}F'_{\theta_0}(x)(2e^{-x}-1/2)dx\right]^2
\end{eqnarray*}
%where $F'_{\theta_0}=\frac{\partial F_{\theta}}{\partial\theta}|_{\theta=\theta_0}$.
We consider three families of Weibull, Linear failure rate and Makeham distributions with the following density functions.
\begin{description}
 \item (1) Weibull Distribution:
\begin{eqnarray*}
f_{\theta}(x)=\theta x^{\theta-1}e^{-x^{\theta}}, \ \ x>0, \
\theta\ge1.
\end{eqnarray*}
%===== ========================================================  LFR
\item (2) Linear Failure Rate Distribution:
\begin{eqnarray*}
f_{\theta}(x)=(1+\theta x)e^{-x-\frac{\theta x^2}{2}}, \ \ x>0, \
\theta\ge0.
\end{eqnarray*}
%=======================================================================  Makeham
 \item (3) Makeham Distribution:
\begin{eqnarray*}
f_{\theta}(x)=(1+\theta (1-e^{-x}))e^{-x-\theta(x+e^{-x}-1)}, \ \
x>0, \ \theta\ge0.
\end{eqnarray*}
\end{description}
  PAE of our test ($\hat{\delta}(F, E_\lambda)$), Deshpande's test ($J_b$), Kochar's test ($T_n$), Link's test ($\Gamma$) and Ahmad's test ($\hat{\Delta}_F$) are presented in Table \ref{PAE}  . In Table \ref{PARE},  PARE  of our test with respect to the others has been obtained. It is observed that our test dominated the others except  Kochar's test for the LFR alternative case.
%($\hat{J}(F_n, 0.4, 3)$ and $\hat{J}(F_n, 0.4, 3)$). 

\begin{table}[!ht]
\caption{ PAE of $\hat{\delta}(F, E_\lambda)$, $J_{0.9}$, $T_n$, $\Gamma$, $\hat\Delta_F$.}\label{PAE}
 % \hrule
\centering
\begin{tabular}{llclllcll}
 Test $\diagdown$ $H_1$& &&Weibull& LFR& Makeham&\\
 \cline{1-2}\cline{3-7}
 $\hat{\delta}(F, E_\lambda)$ &&&1.4414&0.75&0.0833\\
  $J_{0.9}$       &&& 1.35&0.3369&0.0666 \\
  $T_n$          &&&1.247&0.8933&0.0784\\
  $\Gamma$ &&& 1.3867 &0.2681&0.0563\\
  $\hat{\Delta}_F$       &&& 1.35&0.3375&0.0667\\
\hline
\end{tabular}
\end{table}

%======================================================== Table PAER

\begin{table}[!ht] \caption{PARE($\hat{\delta}(F, E), T)=\frac{PAE(\hat{\delta}(F, E_\lambda))}{PAE(T)}$; $ T= J_{0.9}$, $T_n$, $\Gamma$, $\hat\Delta_F$.}\label{PARE}
\vspace{.6cm}
\centering
\begin{tabular}{llccccccc}
 Test$\setminus$ $H_1$& &&Weibull& LFR& Makeham&&\\
 \cline{1-1}\cline{3-8}
  $J_{0.9}$           &&& 1.0677&2.222&1.2489 \\
  $T_n$                &&&1.1558&0.8396&1.0625\\
  $\Gamma$          &&& 1.0394 &2.7974&1.479\\
 $\hat{\Delta}_F$   &&& 1.0677&2.222&1.2489\\
 % $\hat{J}(F_n; 0.4, 3)$   &&&  &1.1510 & 1.0246 \\
 % $\hat{J}(F_n; 0.4, 4)$ &&&  &0.9034 &1.0297\\
\hline
\end{tabular}
\end{table}

In practice, the available samples are small. So, it is important to investigate the  power of  the tests and compare them for small sample sizes.
Proportion of 10000 samples (with small sizes 5(3)15) that  reject exponentiality in favor of IFRA-ness is considered for estimating the power of the tests. In the alternative, we consider Weibull, LFR and Makeham distributions. The critical points of $J_{0.9}$, $T_n$, $\Gamma$ and $\hat{\Delta}_F$ at significance level $\alpha=0.05$ for small sample sizes have been derived from their  corresponding papers.
Table \ref{Tablepowerone} shows the simulated powers of the tests for different alternatives. It is observed that in Weibull and Makeham alternatives our new test is more powerful than the others in all sample sizes.  In the LFR alternative,  Kochar's test  dominates the other tests while  our proposed test is comparable. Also, Kochar's test and Link's test are comparable in Weibull and Makeham alternatives and are more powerful than  the tests of Deshpande and Ahmad.

\begin{table}[htb] 
\centering 
\small 
\caption{ Simulated power of the tests at level of significance 0.05 for  small sample sizes.}\label{Tablepowerone} 
\vspace{0.25cm}
\small\addtolength{\tabcolsep}{-4pt}

\begin{tabular}{cllcl|lll|ccc}
%\hrule
%\hline 
n & &&  {\small Weibull($\theta$)} &  &&  LFR($\theta$) & &&Makeham($\theta$)& \\ 
%& &&  $\theta$ &  &&  LFR ($\theta$) & &&Makeham ($\theta$)&\\
  && 1.2 &2&3 & 0.2 &1&2.5&0.2&1&2.5\\
  \hline         
5 & $\hat{\delta}(F,E_\lambda)$  &0.0902   &0.3886  &0.7496   &0.0645  &0.1029  &0.1422   &0.0579&0.0821&0.1111\\           
   &$J_{0.9}$               &0.0658   &0.1782  &0.3157     &0.0601   &0.0757  &0.0944    &0.0520&0.0685&0.0851 \\ 
   &$T_n$                    &0.0897 &0.3677 &0.7290   &0.0645 &0.1027  &0.1410      &0.0583&0.0833&0.1112\\
   &$\Gamma$              &0.0826 &0.3618  &0.7118   &0.0585   &0.0961& 0.1333  &0.0541&0.0766&0.1037 \\
   &$\hat{\Delta}_F$      &0.0781 &0.3258  &0.6809   &0.059     &0.0637&0.0358   &0.0509&0.0519&0.0391\\  
7 & $\hat{\delta}(F,E_\lambda)$  &0.1117   &0.556   &0.9309      &0.0686   &0.1292  & 0.1854   &0.0617&0.0927&0.1366\\  
   &$J_{0.9}$               &0.0784   &0.2153 &0.4461      &0.0549   &0.0812  &0.1032    &0.0529&0.0666&0.0881\\
   &$T_n$                     &0.1084  &0.5255  &0.9181      &0.0664   &0.1266 &0.1846    &0.0608&0.0906&0.1314\\
   &$\Gamma$              &0.1100     &0.5449  &0.9176      &0.0706   &0.1260   &0.1796      &0.063&0.0974&0.1395 \\
   &$\hat{\Delta}_F$      &0.0988   &0.4600     &0.8758      &0.0602   &0.0806 & 0.0668   &0.0533&0.0612&0.0556\\
9 & $\hat{\delta}(F,E_\lambda)$  &0.1294   &0.7008  &0.9818     &0.0732   &0.1497  &0.2131    &0.0645&0.1087&0.1562\\  
   &$J_{0.9}$               &0.0857   &0.311  &0.6371       &0.0622   &0.0895 &0.1222     &0.0569&0.0797&0.1043\\ 
   &$T_n$                    &0.1227   &0.6705  &0.9763      &0.0741  &0.1469  &0.2087    &0.0658&0.1064&0.1494\\
   &$\Gamma$              &0.1284  &0.6798  &0.9708      &0.0697  &0.1364  &0.1914    &0.0679&0.1072&0.1520 \\
   &$\hat{\Delta}_F$      &0.1068  &0.5828  &0.9476      &0.0630    &0.0952  &0.1032    &0.0557&0.0710&0.0804\\
11 & $\hat{\delta}(F,E_\lambda)$ &0.1440   &0.8038  &0.9967      &0.0767  &0.1773  &0.2872    &0.066&0.1285&0.1994\\ 
   &$J_{0.9}$                &0.0864 &0.3779   &0.7468     &0.0632  &0.1094  &0.1447    &0.0605&0.0795&0.1154\\
   &$T_n$                     &0.1345  &0.7704  &0.9953      &0.0770    &0.1774  &0.2830     &0.0638&0.1227&0.1910\\
   &$\Gamma$              &0.1419  &0.7775  &0.9931      &0.0717   &0.1563  &0.2494    &0.0657&0.1202&0.1844\\
   &$\hat{\Delta}_F$      &0.1204 &0.6867  &0.9857      &0.0695   &0.1193  &0.1452    &0.0571&0.0849&0.1033\\
13 & $\hat{\delta}(F,E_\lambda)$ &0.1629 &0.8777  &0.9989      &0.0882  &0.2042   &0.3229   &0.0716&0.1303&0.2333\\ 
   &$J_{0.9}$                &0.1011  &0.447  &0.8511      &0.0737   &0.1209  &0.1619   &0.0589&0.0875&0.1269\\
   &$T_n$                     &0.1526 &0.8512  &0.9985      &0.0876   &0.2026  &0.3213   &0.0718&0.1247&0.2274\\
   &$\Gamma$              &0.1557 &0.8494  &0.997        &0.0787   &0.1767  &0.2771   &0.0676&0.1215&0.2076\\
   &$\hat{\Delta}_F$      &0.1326 &0.7673  &0.9942      &0.0727   &0.1365  &0.1750   &0.0629&0.0903&0.1278\\     
15 & $\hat{\delta}(F,E_\lambda)$ &0.1764 &0.925   &0.9999      &0.0926   &0.2246  &0.3760   &0.0686&0.1493&0.2554\\
   &$J_{0.9}$               &0.1074 &0.5271  &0.9143      &0.0705  &0.1259   &0.1802  &0.0613&0.0895&0.1430\\
   &$T_n$                     &0.1665 &0.9009 &0.9998      &0.0942   &0.2230  &0.3734  &0.0687&0.1422&0.2511\\
   &$\Gamma$              &0.1800    &0.9067  &0.9995      &0.0862   &0.1903  &0.3217   &0.0664&0.1362&0.2345\\
   &$\hat{\Delta}_F$      &0.1439 &0.8261 &0.9984       &0.0767  &0.1466  &0.2085   &0.0606&0.1005&0.1474\\
\end{tabular}                                                                        
%\end{eqnarray} ""   "0.1094"         ""      ""    ""
\end{table}

\subsection{The Two-Sample}
As mentioned in the introduction, Izadi and Khaledi (2012) proposed and studied a test for the two-sample problem based on kernel density estimation for testing $H_0:F=_*G$ against $H_1: F\le_*G \ \& \ F\neq_*G$. Their test statistic is
\begin{eqnarray*}
  \hat\Delta(F,G)
  %&=&\frac1n\sum_{i=1}^{n}X_i\hat f_n(X_i)- \frac1m\sum_{i=1}^{m}Y_i\hat
  %g_m(Y_i)\\
  =\frac1{n^2a_n}\sum_{i=1}^{n}\sum_{j=1}^{n}X_ik(\frac{X_i-X_j}{a_n})- \frac1{m^2b_m}\sum_{i=1}^{m}\sum_{j=1}^{m}Y_ik(\frac{Y_i-Y_j}{b_m})
\end{eqnarray*}
where $k$ is a known symmetric and bounded density function and $a_n$ and $b_m$ are two sequences of positive real numbers. $k$ and $a_n$ are known as kernel and bandwidth, respectively.

In this section, we compare the empirical power of our new test with the Izadi and Khaledi's test when the kernel, $k$, is the density function of the standard normal distribution and $a_n=n^{-2/5}$ and $b_m=m^{-2/5}$ . We know that the gamma and Weibull family are decreasing  with respect to the shape parameter in the more IFRA order (cf. Marshal and Olkin, 2007, Chapter 9). Also, Izadi and Khaledi (2012) showed that the beta family with density function
\begin{eqnarray}
f(x)=\frac{x^{a-1}(1-x)^{b-1}}{\beta(a,b)}, \ \ x\in[0,1], \ \ a, b>0.
\end{eqnarray}
is increasing with respect to $b$ in the more IFRA order.  So, to evaluate the power of the tests we use the gamma, Weibull and beta families denoted by $G(\alpha,\beta)$, $W(\alpha, \beta)$ and $B(a,b)$, respectively, in the alternative hypothesis. In Table \ref{T3-2}, we generated 10000 samples with sizes $n=m=20, 30, 40, 50, 100$ from distribution $F$ and $G$ given in the table.
We observe that the empirical power of our new test is greater than the empirical power of Izadi and Khaledi's test when $F$ and $G$ belong to Weibull family and is smaller when $F$ and $G$ belong to the gamma and beta families. So, our new test is comparable to Izadi and Khaledi's test.

\begin{table}
\caption{The empirical power of the tests $\hat{\delta}(F,G)$ and $\hat{\Delta}(F, G)$}\label{T3-2}
\begin{center}
  \hrule
\begin{tabular}{cclcccclll}
 &Distribution&&&&&$n=m$&&\\
 \cline{1-2}\cline{5-9}
   $F$&$G$&test&&20&30&40&50&100\\
\hline
    $G(3,1)$&$G(1.5,1)$&$\hat\delta(F,G)$   &&$0.415$ &$0.582$&$0.6764$&$0.7692$&$0.9589$\\
                             &&$\hat{\Delta}(F,G)$   &&$0.470$&$0.642$&$0.7554$&$0.827$&$0.9712$\\                                                                                                         $G(4,1)$&$G(2,1)$&$\hat\delta(F,G)$       &&$0.435$&$0.570$&$0.681$&$0.785$&$0.966$\\
                              &&$\hat{\Delta}(F,G)$  &&$0.492$&$0.672$&$0.755$&$0.825$&$0.968$\\
    $W(3,1)$&$W(1.5,1)$&$\hat\delta(F,G)$  &&0.7946  &0.9346&0.9804&0.9940&$1$\\
                              &&$\hat{\Delta}(F,G)$ &&0.6714   &0.8562&0.9446&0.9796&$1$\\
    $W(4,1)$&$W(2,1)$&$\hat\delta(F,G)$    &&$0.811$ &0.9312&0.9766&$0.993$&$1$\\
                           &&$\hat{\Delta}(F,G)$    &&$0.72$&$0.893$&$0.954$&$0.985$&$1$\\
    $B(1,1.5)$&$B(1,3)$&$\hat\delta(F,G)$   &&$0.1292$&$0.1736$&$0.2332$&$0.2692$&$0.4396$\\
      &&$\hat{\Delta}(F,G)$                         &&$0.1452$&$0.2048$&$0.2806$&$0.3314$&$0.5504$\\
   $B(1.5,2)$&$B(1.5,5)$&$\hat\delta(F,G)$ &&$0.166$&$0.209$&$0.278$&$0.364$&$0.597$\\
      &&$\hat{\Delta}(F,G)$                         &&$0.383$&$0.517$&$0.585$&$0.653$&$0.938$\\
   \hline
\end{tabular}
\end{center}
\end{table}

\section{An application}
In this section we apply our test on a data set from Nelson (1982, page 529) which is a life test to compare two different (old and new) snubber designs. Let $F$ ($G$) be the distribution of lifetime old (new) design population. In Fig. 5, Izadi and Khaledi (2012) plotted TTT-plots  for both  data sets of old and new design. The graphs anticipated IFRA populations for both populations.

Now we apply our IFRA test on the two data sets. Using our one sample test, we get that $\sqrt{12n}(\hat\delta_F-3/2)=-3.579222$ and $\sqrt{12m}(\hat\delta_G-3/2)=-3.085525$ which are less than $-2.376441$ (the critical value at level of significance $\alpha=0.01$ from Table \ref{TabCritic}). So, our test reject exponentiality of both population in favor of IFRA-ness. To compare two populations with respect to more IFRA order, the test statistic  value of the two sample problem is $\sqrt{N}\hat\delta(F,G)/\hat{\sigma}_{F,G}=-0.3762\nless -2.326348=z_{0.01}$. So, at level of significance $\alpha=0.01$, the equality of two populations in more IFRA order is not rejected.

%======================================================= Section coclusion

\section{Summary and Conclusion} 
In order to evaluate the performance of an appliance, we need to compare its aging behavior with some distributions such as exponential, Weibull, gamma,  linear failure rate distributions. The notion of the {\it star order} ( denoted by $\le_*$) is one of the useful tools for this comparison between two distributions. 

In this paper, we  have introduced a new simple test for the problem of testing $H_0: F=_*G$ against $H_1:F\le_*G$ and $F\neq_*G$. 

 In the one-sample problem, let $X_1, \ldots, X_n$ be a random sample from F and  $G=G_0$ where $G_0$ is a known distribution. $H_0$ is rejected at level of significance $\alpha$, for large sample size,  if $n^{1/2}(\hat{\delta}_F-\delta_{G_0})/\sigma_{G_0}< z_{\alpha}$, where
 \begin{eqnarray*}
 \delta_{G_0}=\frac{E_{G_0}[\max\{X_1, X_2\}]}{\mu_{G_0}}, \ \ \  \hat\delta_F=\frac{\displaystyle{\mathop{\sum\sum}_{i\neq j}}\max\{X_i, X_j\}}{n(n-1)\bar X} 
\end{eqnarray*}
and 
\begin{eqnarray*}
    \sigma_{G_0}^2=\frac4{\mu^2_{G_0}}\times
    Var_{G_0}\left(XG_0(X)+\int_{X}^{\infty}tdG_0(t)-\frac{\delta_{G_0}}{2}X\right).
\end{eqnarray*}
In particular, when $G_0$ is an exponential distribution, the null hypothesis in favor of IFRA-ness is  rejected, if $\sqrt{12n}(\hat{\delta}_F-3/2)<z_{\alpha}$.   The exact null distribution of the test statistic has been  obtained and, for small sample sizes 2(1)40, the exact critical points of the test statistics have been computed. Based on Pitman's asymptotic relative efficiency and simulated power, we have compared our test with the tests given by Deshpande (1983), Kochar (1985), Link (1989) and Ahmad (2000).  The results showed that our test relatively dominates the other tests.

In two-sample problem, let $X_1, \ldots, X_n$ and $Y_1, \ldots, Y_m$ be two random samples from $F$ and $G$ respectively. For large sample sizes, we reject $H_0$ in favor of $H_1$ if  
$$\sqrt N(\hat\delta_F-\hat\delta_G)/\hat\sigma_{F,G}<
z_{\alpha}$$
where $N=n+m$ and $\hat\sigma_{F,G}$ has been given in (\ref{SigmaFG}). Using simulation study, we have shown that our test in this case is comparable with the test of Izadi and Khaledi (2012).

\vspace{0.5cm}
\medskip
\noindent{\bf Acknowledgement:}
The authors would like to thank the anonymous associate editor and the referee for
their valuable comments leading to the improvement of our manuscript. The research of  Baha-Eldin Khaledi is partially  supported from Ordered and Spatial Data Center of Excellence of Ferdowsi University of Mashhad.

%============================================================================ References
%\newpage
\vspace{2cm}
  \center {\bf References}
\begin{description}
%\item Ahmad IA (1994) A class of statistics useful in testing
%in increasing failure rate average and new better than used life
%distributions. J Statist Plann Inference 41:141--149

\item Ahmad IA (2000) Testing exponentiality against positive
ageing using kernel methods. \emph{Sankhy\={a} Series A} 62:244--257.

\item Barlow RE,  Proschan F (1981) Statistical theory of
reliability and life testing.\emph{ To Begin With, Silver Spring.}

%\item  Belzunce F,   Candel J and  Ruiz JM  (1998) Testing the stochastic order and the IFR, DFR, NBU, NWU ageing classes ,
%IEEE Trans. Reliab. 47:285–296

\item Deshpande JV (1983) A class of tests for exponentiality against increasing failure rate average alternatives.
\emph{Biometrika} 70:514--518.

\item El-Bassiouny AH (2003) On testing exponentiality against IFRA alternatives. \emph{Appl Math Comput} 146:445--453.

\item Gut A (2009) An Intermediate Course in Probability. \emph{New York: Springer. }

\item Izadi M,  Khaledi BE (2012) On testing more IFRA ordering. \emph{J Statist Plann Inference} 142:840--847.

\item Kochar SC (1985) Testing exponentiality against monotone failure rate average. \emph{Commun Statist Theor Meth} 14:381--392.

\item Kochar SC, Xu M (2011) The tail behavior of the convolutions of Gamma random variables. \emph{J Statist Plann Inference} 141:418--428.

\item Langenberg P, Srinivasan R (1979) Null distribution of the Hollander-Proschan statistic for decreasing mean residual life. \emph{Biometrika} 66:679-980.

\item Lehman EL (1999) Elements of large-sample theory. \emph{New York: Springer-Verlag. }

\item Link WA (1989) Testing for exponentiality against monotone failure rate average alternatives. \emph{Commun Statist Theor Meth} 18:3009--3017.
\item Marshal AW, Olkin I (2007) Life distributions.  \emph{New York: Springer.}

\item Nikitin Y (1995) Asymptotic efficiency of nonparametric tests. \emph{ Cambridge University Press.}

\item Shaked M,   Shanthikumar JG  (2007) Stochastic orders. \emph{New York: Springer-Verlag.}

 \item Xie M,  Lai CD (1996) On the increase of the expected lifetime by parallel redundancy Asia Pac J Oper Res 13:171-179.

\end{description}

\end{document}